\documentclass[12pt]{amsart}
\usepackage{amsmath}
\usepackage{geometry,amsfonts,amssymb,amsthm,txfonts,pxfonts,amscd} 
\def\struckint{\mathop{%
\def\mathpalette##1##2{\mathchoice{##1\displaystyle##2}%
  {##1\textstyle##2}{##1\scriptstyle##2}{##1\scriptscriptstyle##2}}%
\mathpalette
{\vbox\bgroup\baselineskip0pt\lineskiplimit-1000pt\lineskip-1000pt
\halign\bgroup\hfill$}
{##$\hfill\cr{\intop}\cr\diagup\cr\egroup\egroup}%
}\limits}

\newtheorem{theorem}{Theorem}
\newtheorem{lemma}[theorem]{Lemma}
\newtheorem{corollary}[theorem]{Corollary}

\newtheorem{definition}[theorem]{Definition}
\newtheorem{fact}[theorem]{Fact}

\newtheorem{theorem-definition}[theorem]{Theorem-Definition}

\theoremstyle{remark}

\newtheorem{question}[theorem]{Question}

\newcommand{\cx}{\mathbb{C}}

\newcommand{\reals}{\mathbb{R}}
\newcommand{\euc}{\mathbb{E}}

\begin{document}

\title{Simplices and spectra of graphs}
\author{Igor Rivin}
\address{Department of Mathematics, Temple University, Philadelphia}
\email{rivin@temple.edu}
\thanks{The author would like to thank the American Institute of Mathematics for an invitation to the workshop on ``Rigidity and Polyhedral Combinatorics", where this work was started. The author has profited from discussions with Igor Pak, Ezra Miller, and Bob Connelly.}

\date{\today} 


\begin{abstract}
In this note we show the $n-2$-dimensional volumes of codimension $2$ faces of an $n$-dimensional simplex are algebraically independent functions of the lengths of edges. In order to prove this we compute the complete spectrum of a combinatorially interesting graph.
\end{abstract}

\maketitle

\section*{Introduction}

Let $T_n$ be the set of congruence classes of simplices in Euclidean space $\euc^n.$ The set $T_n$ is an open manifold (also a semi-algebraic set) of dimension $(n+1)n/2.$ Coincidentally, a simplex $T\in T_n$ is determined by the $(n+1)n/2$ lengths of its edges. Furthermore, the square of the volume of $T \in T_n$ is a polynomial in the squares of the edgelengths. This polynomial is given by the Cayley-Menger determinant formula:
\[
V^2(T) = \dfrac{(-1)^{n+1}}{2^n (j!)^n} \det C,
\]
where $C$ is the \emph{Cayley-Menger matrix,} defined as follows:

\[
C_{ij} = \begin{cases}
0,& i=j\\
1, & \mbox{if $i=1$ or $j=1$, but not both}\\
l_{(i-1)(j-1)}^2,& \mbox{otherwise}
\end{cases}
\]

Note that an $n$-dimensional simplex also has $(n+1)n$ $n-2$ dimensional faces, and so the following question is natural:

\begin{question}
\label{warrenq}
Is the congruence class of the simplex $T$ determined by the $n-2$-dimensional volumes of the $n-2$-dimensional faces?
\end{question}

Question \ref{warrenq} must be classical, but the first reference that I am aware of is Warren Smith's Princeton PhD thesis \cite{warrenthesis}.

In the AIM workshop on Rigidity and Polyhedral Combinatorics Bob Connelly (who was unaware of the reference \cite{warrenthesis}) raised the following:

\begin{question}
\label{connellyq}
Is the \emph{volume} of the simplex $T$ determined by the $n-2$-dimensional volumes of the $n-2$-dimensional faces?
\end{question}
In fact, Connelly stated Question \ref{connellyq} for $n=4,$ which is the first case where the question is open (for $n=3$ the answer is trivially "Yes", since $3-2=1,$ and we are simply asking if the volume of the simplex is determined by its edge-lengths. In dimension 2, the answer is trivially "No", since $2-2=0,$ and the volume of codimension-2 faces of a triangle carries no information.

Clearly, the affirmative answer to Question \ref{warrenq} implies an affirmative answer to Question \ref{connellyq}. At the time of this writing, both questions are open. In this note we show

\begin{theorem}
\label{mainthm}
The $(n+1)n/2$ $n-2$-dimensional volumes of the $n-2$-dimensional faces of an $n$-dimensional simplex are algebraically independent over $\cx[l_{12}, \dotsc, l_{((n+1)(n+1}].$
\end{theorem}
Theorem \ref{mainthm} is clearly a necessary step in the direction of resolving Question \ref{warrenq}, but is far from sufficient. To show it, consider the map of
$\reals^{(n+!)n/2}$ to $\reals^{(n+1)n/2},$ which sends the vector $E$ of edge-lengths to the vector $F$ of areas of $n-2$-dimensional faces. To show Theorem \ref{mainthm} it is enough to check that the Jacobian $J(E)=\partial F/\partial E$ is non-singular at \emph{one} point. We will use the most obvious point $p_1:$ one corresponding to a regular simplex with all edge-lengths equal to $1.$  By symmetry considerations, the Jacobian $J(p_1)$ can be written as
\[
J(p_1) = c M,
\]
where 
\[
M_ij = 
\begin{cases}
1,& \mbox{if the $j$-th edge is incident to the $i$-th face of dimension $n-2.$}\\
0,& \mbox{otherwise}
\end{cases}
\]
The first observation is
\begin{lemma}
The constant $c$ above is not equal to $0.$
\end{lemma}
\begin{proof}
This follows from the observation that the volume of a $k$-dimensional simplex is a \emph{homogeneous} function of the edge-lengths, of degree $k.$ An application of Euler's Homogeneous Function Theorem shows that at $p_1,$ 
\[
\dfrac{\partial F_i}{\partial e_j} = \begin{cases} \frac2{n-1}F,& \mbox{if $e_j$ is incident to $F_i$}\\
0,& \mbox{otherwise},
\end{cases}
\]
where $F$ is the common value of the $n-2$-dimensional volume of the $F_i,$ which implies that $c=\frac{2}{n-1} F.$
\end{proof}

Theorem \ref{mainthm} thus reduces to the assertion that the determinant of the matrix $M$ is not zero. We will be able to prove a much stronger result (of interest in its own right):

\begin{theorem}
\label{singvals}
The singular values of $M$ are as follows: 
The value $1$ appears with multiplicity $(n+1)(n-1)/2.$ The value $(n-2)$ appears with multiplicity $n.$ The value $(n-2)(n-1)$ appears once.
\end{theorem}

\begin{corollary}
\label{thedet}
The absolute value of the determinant of $M$ equals
\[
(n-2)^{n+1}(n-1) \neq 0,
\]
for $n>2.$
\end{corollary}

To prove Theorem \ref{singvals}, first recall that the singular values of $M$ are the positive square roots of the eigenvalues of $N=MM^t.$  In its turn, $N$ has rows and columns indexed by $n-1$-subsets of $R_{n+1} = \{1, \dotsc, n+1\}.$ The $ij$-th element of $N$ equals the number of $2$-element subsets the $i$-th and the $j$-th two element subsets $s_i$ and $s_j$ have in common. This, in turn, can be written as follows:
\[
N_{ij} = \begin{cases}n(n-1)/2,& i=j\\
	                           (n-1)(n-2)/2,& |s_i\cap s_j|= n-2\\
	                           (n-2)(n-3)/2, & |s_i\cap s_j| = n-3.
	                           \end{cases}
	                           \]
	                           The matrix $N$ is the adjacency matrix of a multi-graph $G_N$, which has a rather large symmetry group. These symmetries will allow us to obtain the complete spectral decomposition of the matrix $N$. Indeed, the symmetry group of $G_N$ is the symmetric group $S_{n+1},$ while the stabilizer $\Gamma_i$ of a vertex $s_i$ is the group $S_{n-1}\times S_2$  (the first factor acts on $s_i$ itself, the second on $R_{n+1} \backslash s_i.$) The action of $\Gamma_i$ on $G_N$ has three orbits. The first consists of $s_i$ itself. The second consists of all $s_j$ such that $|s_i\cap s_j| = n-2,$ the third of all $s_j$ such that $|s_i\cap s_j| = n-3.$
	                           
	                           At this point it behooves us to recall the concept of \emph{graph divisor}.


\section{Graph divisors}
The concept of graph divisor is discussed at great length in the books \cite{cvetkovicspectra} and \cite{cvetkoviceigen}
\begin{definition}
Given an $s\times s$ matrix $B = (b_{ij}),$ let the vertex set of a (multi)graph $G$ be partitioned into (non-empty) subsets $X_1, X_2, \dotsc, X_s,$ so that for any $i, j = 1, 2, \dotsc, s,$ each vertex from $X_i$ is adjacent to exactly $b_{ij}$ vertices of $X_j.$ The multidigraph $H$ with adjacency matrix $B$ is called a \emph{front divisor} of $G.$
\end{definition}

The importance of this concept to our needs is that the characteristic polynomial of a graph divisor \emph{divides} the characteristic polynomial of the adjacency matrix of $G.$ (hence the name). The most interesting (to us, anyway) example of a graph divisor arises by having a subgroup $\Gamma$ of the automorphism group of $G.$ The quotient of $G$ by $\Gamma$ is the divisor we consider. Every eigenvector of $\Gamma \backslash G$ lifts to an eigvenvector of $G$ with the same eigenvalue.

If we consider our graph $G_N$ and the action of $\Gamma_i,$ we observe (after a rather tedious computation) that the front divisor corresponding to the action of $\Gamma_i$ on $G_N$ has adjacency matrix
\[
D=\begin{pmatrix}
 \frac{1}{2} (n-2) (n-1) & (n-3) (n-2) (n-1) & \frac{1}{4} (n-4) (n-3)
   (n-2) (n-1) \\
 \frac{1}{2} (n-3) (n-2) & n^3-7 n^2+17 n-14 & \frac{1}{4} \left(n^4-10
   n^3+39 n^2-70 n+48\right) \\
 \frac{1}{2} (n-4) (n-3) & n^3-8 n^2+23 n-24 & \frac{1}{4} \left(n^4-10
   n^3+43 n^2-90 n+76\right)
\end{pmatrix}
\]
A simple computation shows that the eigenvalues of $D$ are $(n-1)^2(n-2)^2, (n-2)^2, 1,$ while the corresponding eigenvectors are (respectively):

\begin{gather*}
(1, 1, 1),\\ ((1-n)/2, (3-n)/4, 1),\\ (
(2-3n+n^2).2, (2-n)/2, 1).
\end{gather*}

Graph divisors are a very useful tool, but they have too obvious shortcomings:

\begin{enumerate}
\item Not all eigenvalues of the graph $G$ are captured by the graph divisor.
\item We have no information on the \emph{multiplicity} of any of the eigenvalues that are captured.
\end{enumerate}

Here, however, we have a \emph{deus ex machina} in the form of 
\section{Gelfand pairs}
We will not need any more than the (rather little) presented in Diaconis' little book \cite[page 54]{diaconisbook}

First:
\begin{definition}
Let $G$ be a group acting transitively on a finite set $X$ with isotropy group $N.$ A function $f: G\rightarrow \cx$ is called \emph{$N$-bi-invariant} if $f(n_1 g n_2) = f(g),$ for all $n_1, n_2 \in N, g\in G.$ The pair $G, N$ is called a \emph{Gelfand pair} if the convolution of $N$-bi-invariant functions is commutative.
\end{definition}

In our application, $G = S_{n+1},$ $N = S_{n-1}\times S_2,$ and $X=G_N.$

The further results we need are the following (the citations are from \cite{diaconisbook}:
\begin{lemma}[Lemma 5, page 53]
\label{lemma5}
Let $\rho, V$ be an irreducible representation of the finite group $G.$ Let $N \subset G$ be a subgroup and let $X=G/N$ be the associated homogeneous space. The number of times that $\rho$ appears in $L^2(X)$ equals the dimension of the space of $N$-fixed vectors in $\rho, V$.
\end{lemma}

\begin{theorem}[Theorem 9, page 54]
The following three conditions are equivalent:
\begin{enumerate}
\item $G, N$ is a Gelfand pair.
\item The decomposition of $L^2(X)$ into irreducible representations of $G$ is multiplicity-free.
\item For every irreducible representation $(\rho, V)$ there is basis of $V$ such that $\hat{f}(\rho) = \bigl( \begin{smallmatrix}
* & 0 \\ 0 & 0\end{smallmatrix}\bigr),$ for all $N$-bi-invariant functions $f.$
\end{enumerate}
\end{theorem}

The main significance (as seen by this author) of the Gelfand pair property is that the eigenvalues of $N\backslash G / N$ (which is the "front divisor" $D$) are the same as the eigenvalues of  $X= G/N,$ since each irreducible factor of $L^2(X)$ contains an invariant function (as per Lemma \ref{lemma5})

We will finally need the following (attributed by Diaconis to \cite{jamessn}):
\begin{fact}
\label{jamesfact}
The pair $S_{n+1}, S_{n-1}\times S_2$ is a Gelfand pair. The dimensions of the irreducible representations of $G_n$ are $(n+1)n/2, n, 1.$
\end{fact}

We now have everything we need to finish the proof of Theorem \ref{singvals}, and hence Theorem \ref{mainthm}. Since the graph laplacian (or adjacency matrix) commutes with the action of the automorphism group, each eigenspace of the adjacency matrix $M$ is a sum of irreducible representations of $S_{n+1}.$ By Fact \ref{jamesfact} there are precisely three eigenspaces, and their dimensions are as stated in Theorem \ref{singvals}. The only question is how to decide which eigenspace has which dimension. Since we know the eigenvectors of the ``front divisor'' matrix $D$ (see the end of the Introduction) we can match them up with the ``spherical functions'' of the Gelfand pair $S_{n+1}, S_{n-1}\times S_2$ (as given on \cite[page 57]{diaconisbook}, whence the result follows.

\bibliographystyle{plain}
\bibliography{rivin}
\end{document}